\documentclass[11pt]{amsart}

\usepackage{pstricks,graphicx,pst-node}
\newtheorem{theorem}{Theorem}[section]

\newtheorem{proposition}[theorem]{Proposition}
\newtheorem{corollary}[theorem]{Corollary}

\newtheorem{remark}[theorem]{Remark}

\newenvironment{rmkT}[1]{\noindent\textit{Remarks on \reftab{#1}.}}{\vspace{\baselineskip}}

\newcommand{\reftab}[1]{Table \ref{tab:#1}}
\newcommand{\refeq}[1]{(\ref{eqn:#1})}

\newcommand{\refthm}[1]{Theorem \ref{thm:#1}}

\newcommand{\refprop}[1]{Proposition \ref{prop:#1}}

\newcommand{\refsec}[1]{Section \ref{sec:#1}}

\newcommand{\Sn}{\mathcal{S}_n}

\newcommand{\head}{\mathrm{h}}
\newcommand{\tail}{\mathrm{t}}

\newcommand{\SEP}{\mathrm{SEP}}

\newlength{\hsp}
\newlength{\vsp}
\newlength{\vspi}
\begin{document}

%%%%%%%%%%%%%%%%%%%%%%%%%%%%%%%%%%%%%%%%%%%%%%%%%%%%%%%%%%%%
%  TITLE PAGE information
%%%%%%%%%%%%%%%%%%%%%%%%%%%%%%%%%%%%%%%%%%%%%%%%%%%%%%%%%%%%

%     [Short Title]{Full Title}
\title[Riffle shuffles with repeated cards]{Riffle shuffles of a deck with repeated cards}

%    Information for first author
\author[Assaf]{Sami Assaf}
\address{Department of Mathematics, Massachusetts Institute of Technology, 77 Massachusetts Avenue, Cambridge, MA 02139-4307}
%\curraddr{}
\email{sassaf@math.mit.edu}
%\thanks{}

%    Information for second author
\author[Diaconis]{Persi Diaconis}
\address{Department of Statistics, Stanford University, 390 Serra Mall, Stanford, CA 94305-4065}
%\curraddr{ }
%\email{}
%\thanks{}

%    Information for third author
\author[Soundararajan]{K. Soundararajan}
\address{Department of Mathematics, Stanford University, 450 Serra Mall, Building 380, Stanford, CA 94305-2125}
%\curraddr{}
\email{ksound@math.stanford.edu}
%\thanks{}

%    General info
%\subjclass[2000]{Primary ; Secondary }

\date{\today}

%\dedicatory{}

%\keywords{riffle shuffles}

\begin{abstract}
  We study the Gilbert-Shannon-Reeds model for riffle shuffles and ask
  'How many times must a deck of cards be shuffled for the deck to be
  in close to random order?'. In 1992, Bayer and Diaconis gave a
  solution which gives exact and asymptotic results for all decks of
  practical interest, e.g. a deck of 52 cards. But what if one only
  cares about the colors of the cards or disregards the suits focusing
  solely on the ranks? More generally, how does the rate of
  convergence of a Markov chain change if we are interested in only
  certain features? Our exploration of this problem takes us through
  random walks on groups and their cosets, discovering along the way
  exact formulas leading to interesting combinatorics, an 'amazing
  matrix', and new analytic methods which produce a completely general
  asymptotic solution that is remarkable accurate.
\end{abstract}

\maketitle

%%%%%%%%%%%%%%%%%%%%%%%%%%%%%%%%%%%%%%%%%%%%%%%%%%%%%%%%%%%%
\section{Introduction}
%%%%%%%%%%%%%%%%%%%%%%%%%%%%%%%%%%%%%%%%%%%%%%%%%%%%%%%%%%%%
\label{sec:intro}

A basic question in scientific computing is `How many times must an
iterative procedure be run?'. A basic answer is `It depends.'. In this
paper we study the mixing properties of the Gilbert-Shannon-Reeds
model \cite{Gilbert1955,Reeds1976} for riffle shuffling a deck of $n$
cards and ask how many times the deck must be shuffled for the cards
to be in close to random order. Our answer depends not only on the
metric we use to measure distance to uniformity, but also on the
particular properties of the deck that are of interest. 

To be precise, we consider a `deck' to be a multiset of $n$ cards. We
shuffle the deck by first cutting it into two piles according to the
binomial distribution, and then riffling the piles together by
successively dropping cards from either pile with probability
proportional to the size. This process defines a measure, denoted
$Q_2(\sigma)$, on the symmetric group $\Sn$. Repeated shuffles are
defined by {\em convolution powers}
\begin{equation}
  Q_2^{*k}(\sigma) = \sum_{\omega \cdot \tau = \sigma} Q_2(\tau)
  Q_2^{*(k-1)}(\omega) .
\label{eqn:convolution}
\end{equation}

This shuffling model, which accurately models how most people actually
shuffle a deck of cards, was introduced by Gilbert and Shannon
\cite{Gilbert1955} and independently by Reeds \cite{Reeds1976}.

Bayer and Diaconis \cite{BaDi1992} generalized this to $a$-shuffles,
which is the natural extension to shuffling with $a$ hands: the deck
is cut into $a$ packets by multinomial distribution and cards are
successively dropped from packets with probability proportional to
packet size. Letting $Q_a(\sigma)$ denote this measure, they show that
convolution of general $a$-shuffles is as nice as possible, namely
\begin{equation}
  Q_a * Q_b = Q_{ab}.
\label{eqn:conv}
\end{equation}
Thus it is enough to study a single $a$-shuffle of the deck.

To that end, denote the \textit{uniform distribution} by $U =
U(\sigma)$. For a deck with $n$ distinct cards, $U = 1/n!$, and for a
more general deck with $D_1$ $1$'s, $D_2$ $2$'s, up to $D_m$ $m$'s, we
have $U = 1/\binom{D_1 + \cdots + D_m}{D_1, \ldots, D_m}$. There are
several ways to measure the distance between $Q_a$ and $U$, though for
the purposes of this paper we restrict our attention to total
variation distance and separation distance.

The {\em total variation distance} is defined by
\begin{equation}
  \| Q_a - U \|_{TV} = \max_{\mathrm{subsets} A} |Q_a (A) - U(A)| =
  \frac{1}{2} \sum_{\sigma} |Q_a(\sigma) - U(\sigma)| .
\label{eqn:TVdist}
\end{equation}
In general, the formulas for $Q_a(\sigma)$ may be quite complicated,
making calculations of total variation intractable. Therefore we will
also consider the \textit{separation distance} defined by
\begin{equation}
  \SEP(a) = \max_{\sigma} 1 - \frac{Q_a (\sigma)}{U(\sigma)} .
\label{eqn:SEP}
\end{equation}
Here, only a single probability needs to be computed, though as we
shall see even that can be difficult. From the formulas above, one can
easily see that separation provides an upper bound for total
variation, which makes separation a good measure to use when total
variation becomes too complicated to compute.

In widely cited works, Aldous \cite{AlDi1986} and Bayer and Diaconis
\cite{BaDi1992} show that $\frac{3}{2} \log_2(n) + c$ shuffles are
necessary and sufficient to make the total variation distance small,
while $2\log_2(n)+c$ shuffles are necessary and sufficient to make
separation small. These results, however, look at all aspects of a
permutation, i.e. consider a deck with distinct cards. In many card
games, only certain aspects of the permutation matter. For instance,
in Baccarat, suits are irrelevant and all $10$'s and picture cards are
equivalent, and in ESP card guessing experiments, a Zener deck of $25$
cards with each of $5$ symbols repeated five times is used. It is
natural, therefore, to ask how many shuffles are required in these
situations, and so we consider a deck to have repeated cards.

Many results are known for how long it takes certain features of a
permutation, e.g. longest cycle, descent structure, etc, to become
random; for a thorough treatment of such results see
\cite{Diaconis2003}. The particular problem we address in this paper
was first addressed by Conger and Viswanath \cite{CoVi2006,CoVi2007}
who derive remarkable numerical procedures giving useful answers for
cases of practical interest.

In this paper, we present many of our main results from
\cite{ADS2008}, giving exact formulae and asymptotics for a deck of
$n$ cards with $D_1$ cards labelled $1$, $D_2$ cards labelled $2$,
$\ldots$, $D_m$ cards labelled $m$. Our results are proved from the
deck starting `in order', i.e. with $1$'s on top through $m$'s at the
bottom. In \refsec{groups}, we show that the processes we study are
Markov by framing the problem in the context of random walks on
cosets. We derive a formula for the transition matrix following a
single card in \refsec{Aspades}, and show that this matrix shares many
properties with Holte's `Amazing Matrix' \cite{Holte1997}. In
\refsec{deck}, we consider a general deck, limiting our metric to the
separation distance, and derive new formulae and asymptotic
approximations which we unify into our `rule of thumb'
formula. \refsec{initial} shows that our results depend on the initial
configuration of the deck, a fact also observed by Conger and
Viswanath \cite{CoVi2006,CoVi2007,CoVi}. This extended abstract contains
precise statements of our main results along with the main ideas of
the proofs; for full details see \cite{ADS2008}.

%%%%%%%%%%%%%%%%%%%%%%%%%%%%%%%%%%%%%%%%%%%%%%%%%%%%%%%%%%%%
\section{Random walks on Young subgroups}
%%%%%%%%%%%%%%%%%%%%%%%%%%%%%%%%%%%%%%%%%%%%%%%%%%%%%%%%%%%%
\label{sec:groups}

In this section, we reformulate shuffling in terms of random walks on
a finite group, so that our investigation of particular properties of
a deck becomes a quotient walk on Young subgroups of $\Sn$.

Let $G$ be a finite group, and let $Q$ be a probability on $G$, i.e.
$Q(g) \geq 0$ and $\sum_{g \in G} Q(g) = 1$. Take a \textit{random
  walk on $G$} by repeatedly choosing elements independently from $G$
with probability $Q$, say $g_1, g_2, g_3, \ldots$, and, beginning with
the identity element $1_{G}$, multiply on the left by $g_i$. This
generates the following sequence of elements, the left walk,
\begin{displaymath}
  1_G, \ g_1, \ g_2 g_1, \ g_3 g_2 g_1, \ \ldots .
\end{displaymath}
By inspection, the chance that the walk is at $g$ after $k$ steps is
given by convolution formula (\ref{eqn:convolution}) $Q^{*k}(g)$,
where $Q^0(g) = \delta_{1_G,g}$.

To focus on certain aspects of the walk, we choose a subgroup and
consider the \textit{quotient walk} as follows. Let $H \leq G$ be a
subgroup of $G$, and let $X$ denote the set of left cosets of $H$ in
$G$, i.e. $X = G/H = \{xH\}$. The quotient walk on $X$ is derived from
the left walk on $G$ by reporting the coset to which the current
position of the walk belongs. This defines a Markov chain on $X$ with
transition matrix given by
\begin{equation}
  K(x,y) = Q(y H x^{-1}) = \sum_{h \in H} Q(yhx^{-1}) .
\label{eqn:transition}
\end{equation}
Note that $K$ is well-defined (i.e. independent of the choice of coset
representatives) and doubly stochastic. Thus the uniform distribution
on $X$, $U = |H|/|G|$, is a stationary distribution for $K$. The
following result, showing that powers of $K$ correspond precisely to
convolving and taking cosets, is intuitively obvious with a
straightforward proof.

\begin{proposition}
  For $Q$ a probability distribution on a finite group $G$ and $K$ as
  defined in \refeq{transition}, we have
  \begin{displaymath}
    K^{l}(x,y) = Q^{*l} (yHx^{-1}) .
  \end{displaymath}
\label{prop:convolve}
\end{proposition}

We may identify permutations in $\Sn$ with arrangements of a deck of
$n$ cards by setting $\sigma(i)$ to be the label of the card at
position $i$ from the top. For instance, the permutation $2 \ 1 \ 4 \
3$ is associated with four cards where ``2'' is on top, followed by
``1'', followed by ``4'', and finally ``3'' is on the
bottom. Therefore the random walk on $\Sn$ with probability $Q_2$
corresponds precisely to riffle shuffles of a deck of $n$ distinct
cards.  If we consider the first $D_1$ cards to be labelled ``1'', the
next $D_2$ cards to be labelled ``2'', and so on up to the last $D_m$
cards labelled ``$m$'', then this corresponds precisely to the coset
space of a Young subgroup,
\begin{displaymath}
  X = \Sn / \left(\mathcal{S}_{D_1} \times \mathcal{S}_{D_2} \times
    \cdots \times \mathcal{S}_{D_m}\right).
\end{displaymath}
Thus \refprop{convolve} shows that the processes studied in the body
of this paper are Markov chains.

%%%%%%%%%%%%%%%%%%%%%%%%%%%%%%%%%%%%%%%%%%%%%%%%%%%%%%%%%%%%
\section{A new `amazing' matrix}
%%%%%%%%%%%%%%%%%%%%%%%%%%%%%%%%%%%%%%%%%%%%%%%%%%%%%%%%%%%%
\label{sec:Aspades}

Suppose the ace of spades is on the bottom of a deck of $n$ cards. How
many shuffles does it take until this one card is close to uniformly
distributed on $\{1,2,\ldots,n\}$? We analyze this problem by writing
down the transition matrix following a single card through an
otherwise indistinguishable deck.

\begin{proposition}
  Let $P_a(i,j)$ be the chance that the card at position $i$ moves
  to position $j$ after an $a$-shuffle. For $1\leq i,j\leq n$,
  $P_a(i,j)$ is given by
  \begin{displaymath}
    \frac{1}{a^n} \sum_{k=1}^{a} \sum_{r=l}^{u} \!
    \binom{j\!-\!1}{r} \! \binom{n\!-\!j}{i\!-\!r\!-\!1} k^r
    (a-k)^{j-1-r}(k-1)^{i-1-r}(a-k+1)^{(n-j)-(i-r-1)}
  \end{displaymath}
  where $r$ ranges from $l = \max(0,(i+j)-(n+1))$ to $u = \min(i-1,j-1)$.
\label{prop:paij}
\end{proposition}

\begin{proof}
  Consider the number of ways that an inverse $a$-shuffle can bring
  the card at position $j$ to position $i$. First, the card at
  position $j$ must have come from some pile, say $k$, $1 \leq k \leq
  a$. Say $r$ of the cards above this came from piles $1$ to $k$, and
  so the remaining $j-1-r$ came from piles $k+1$ to $a$. Those $r$
  cards all must appear before the card at position $j$ in
  $\binom{j-1}{r}$ ways. This leaves $i-1-r$ cards below position $j$
  which came from piles $1$ to $k-1$ in $\binom{n-j}{i-r-1}$ ways, and
  the remaining cards must be from piles $k$ to $a$.
\end{proof}

For example, the $n\times n$ transition matrices for $n=2,3$ are given
below. 
\begin{displaymath}
  \begin{array}{ccc}
    \displaystyle{\frac{1}{2a} \left( \begin{array}{cc}
        a+1 & a-1 \\ a-1 & a+1
      \end{array} \right)} & \hspace{2em} &
    \displaystyle{\frac{1}{6a^2} \left( \begin{array}{ccc}
        (a+1)(2a+1) & 2(a^2-1) & (a-1)(2a-1) \\
        2(a^2-1)    & 2(a^2+2) & 2(a^2-1) \\
        (a-1)(2a-1) & 2(a^2-1) & (a+1)(2a+1)
      \end{array} \right)}
  \end{array}
\end{displaymath}

These matrices share many properties, given in \refprop{amazing},
with the `amazing matrix' discovered by Holte \cite{Holte1997} in his
study of the `carries process' of ordinary addition. Diaconis and
Fulman \cite{DiFu2008} show that Holte's matrix is also the transition
matrix for the number of descents in repeated $a$-shuffles. We have
not been able to find a closer connection between the two matrices.

\begin{proposition} The transition matrices following a single card
  have the following properties:
  \begin{enumerate}
  \item they are {\em cross-symmetric}, i.e. $P_a(i,j) =
    P_a(n-i+1,n-j+1)$;
  \item they are multiplicative, i.e. $P_a \cdot P_b = P_{ab}$;
  \item the eigenvalues form the geometric series $1, 1/a, 1/a^2,
    \ldots, 1/a^{n-1}$;
  \item the right eigen vectors are independent of $a$ and have the
    simple form: \\
    $V_m(i) = (i-1)^{i-1} \binom{m-1}{i-1} + (-1)^{n-i+m}
    \binom{m-1}{n-i}$ for $1/a^m$, $m \geq 1$.
  \end{enumerate}
\label{prop:amazing}  
\end{proposition}

\begin{proof}
  The cross-symmetry (1) follows from \refprop{paij}, and the
  multiplicative property (2) follows from the shuffling
  interpretation and equation \refeq{conv}. Property (1) implies that
  the eigen structure is quite constrained. Properties (3) and (4)
  follow from results of Cuicu \cite{Ciucu1998}.
\end{proof}

The following Corollary also follows as a special case of Theorem 2.2
in \cite{CoVi2006}.

\begin{corollary}
  Consider a deck of $n$ cards with the ace of spades starting at the
  bottom. The chance that the ace of spades is at position $j$
  from the top after an $a$-shuffle is
  \begin{equation}
    Q_a(j) \ = P_a(n,j) \ = \ \frac{1}{a^n} \sum_{k=1}^{a} (k-1)^{n-j} k^{j-1}.
    \label{eqn:Qai}
  \end{equation}
  \label{cor:Qai}
\end{corollary}

From the explicit formula, we are able to give exact numerical
calculations and sharp asymptotics for any of the distances to
uniformity. The results below show that $\log_2 n +c$ shuffles are
necessary and sufficient for both separation and total variation (and
there is a cutoff for these). This is surprising since, on the full
permutation group, separation requires $2\log_2 n+c$ steps whereas
total variation requires $\frac{3}{2}\log_2 n +c$. Of course, for any
specific $n$, these asymptotic results are just indicative.

\newcommand{\rrc}{@{\extracolsep{10pt}}r}

\begin{table}[ht]
  \begin{center}
    \caption{\label{tab:AS} Distance to uniformity for a deck of $52$
      cards. The upper table assumes distinct cards, and the lower
      table follows a single card starting at the bottom of the deck.}
    \begin{tabular}{c | cccc cccc cccc}
      \hline
      & 1 & 2 & 3 & 4 & 5 & 6 & 7 & 8 & 9 & 10 & 11 & 12 \\ \hline \\
      [-.5\vsp]  
      $TV$ & 1.00 & 1.00 & 1.00 & 1.00 & .924 & .614 & .334 &
      .167 & .085 & .043 & .021 & .010 \\ [.5\vsp]
      $\SEP$ & 1.00 & 1.00 & 1.00 & 1.00 & 1.00 & 1.00 & 1.00 & .996 &
      .931 & .732 & .479 & .278 \\  \hline
    \end{tabular}

    \vspace{\baselineskip}
    \begin{tabular}{c | cccc cccc cccc}
      \hline
       & 1 & 2 & 3 & 4 & 5 & 6 & 7 & 8 & 9 & 10 & 11 & 12 \\ \hline \\
       [-.5\vsp]  
      $TV$ & .873 & .752 & .577 & .367 & .200 & .103 & .052 &
      .026 & .013 & .007 & .003 & .002 \\ [.5\vsp]
      $\SEP$ & 1.00 & 1.00 & .993 & .875 & .605 & .353 & .190 &
      .098 & .050 & .025 & .013 & .006 \\ \hline
    \end{tabular}
%%
%%    \vspace{\baselineskip}
%%    \begin{tabular}{c | cccc cccc cccc}
%%      \hline
%%       & 1 & 2 & 3 & 4 & & & & & & & & \\ \hline \\
%%       [-.5\vsp]  
%%      $TV$ & .494 & .152 & .001 & .000 & \white{.000} & \white{.000} &
%%      \white{.000} & \white{.000} & \white{.000} & \white{.000} &
%%      \white{.000} & \white{.000} \\ [.5\vsp] 
%%      $\SEP$ & 1.00 & .487 & .003 & .000 & \white{.000} & \white{.000}
%%      & \white{.000} & \white{.000} & \white{.000} & \white{.000} &
%%      \white{.000} & \white{.000} \\ \hline 
%%    \end{tabular}
  \end{center}
\end{table}

\begin{rmkT}{AS}
  We use \refprop{paij} to give exact results when $n=52$. For
  comparison, the upper table gives exact results for the full deck
  using \cite{BaDi1992}. The lower table shows that it takes about
  half as many shuffles to achieve a given degree of mixing for a card
  at the bottom of the deck. For example, the widely cited `$7$
  shuffles' for total variation drops this distance to $.334$ for the
  full ordering, but this requires only $4$ shuffles to achieve a
  similar degree of randomness for a single card at the bottom.
\end{rmkT}

For asymptotic results, we first derive an approximation to
separation, which also serves as an upper bound for total
variation. Finally, we derive a matching lower bound for total
variation. Proofs have been omitted for brevity, but again full
details are available in \cite{ADS2008}.

\begin{proposition}
  After an $a$-shuffle, the probability that the bottom card 
  is at position $i$ satisfies 
  \begin{displaymath}
    \frac{1}{a} \frac{\alpha^{n-i+1}}{1-\alpha^n} \le   Q_a(i) \le 
    \frac{1}{a} \frac{\alpha^{n-i}}{1-\alpha^{n-1}}, 
  \end{displaymath}
  where for brevity we have set $\alpha = 1-1/a$.  In particular, the 
  separation distance satisfies 
  \begin{displaymath}
    1- \frac{n}{a} \frac{\alpha^{n}}{1-\alpha^n} \le \SEP(a) \le
    1-\frac{n}{a} \frac{\alpha^{n-1}}{1-\alpha^{n-1}}.
  \end{displaymath}
  \label{prop:sepi-bound}
\end{proposition}

If $a = 2^{\log_2(n)+c} = n2^c$, then our result shows that the
$\SEP(a)$ is approximately
\begin{displaymath}
  1 -\frac{1}{2^c} \frac{e^{-2^{-c}}}{1-e^{-2^{-c}}}, 
\end{displaymath}
and for large $c$ this is $\approx 2^{-c-1}$.  The fit to the data in
\reftab{AS} is excellent: for example after ten shuffles of a
fifty-two card deck we have $2^{-c-1} = \frac{26}{1024}$ which is very
nearly the observed separation distance of $0.025$.

\begin{remark}
  \refprop{sepi-bound} gives a local limit for the probability that
  the original bottom card is at position $j$ \textit{from the
    bottom}. When the number of shuffles is $\log_2 n +c$, the density
  of this (with respect to the uniform measure) is asymptotically
  $z(c) e^{-j/2^c}$, with $z$ a normalizing constant ($z(c) =
  1/2^c(e^{j/2^c}-1)$). The result is uniform in $j$ for $c$ fixed,
  $n$ large.
\end{remark}

\begin{proposition} Consider a deck of $n$ cards with the ace of
  spades at the bottom.  With $\alpha= 1-1/a$, the total variation
  distance for the mixing of the ace of spades after an $a$-shuffle is
  at most
  \begin{displaymath}
    \frac{\alpha^{n+1}}{1-\alpha^n}
    -\frac{a\alpha^2(1-\alpha^{n-1})}{n(1-\alpha^n)}  +\frac{1}{n\log
      (1/\alpha)} \log \left(\frac an
      \frac{1-\alpha^n}{\alpha^{n+1}}\right), 
  \end{displaymath}
  and at least 
  \[ 
  \frac{\alpha^n}{1-\alpha^{n-1}} - \frac{a(1-\alpha^n)}{n\alpha(1-\alpha^{n-1})} 
  + \frac{1}{n\log (1/\alpha)} \log \Big( \frac{a}{n} \frac{1-\alpha^{n-1}}{\alpha^{n-1}}\Big).
  \]
\label{prop:TV-bound}
\end{proposition}

After $\log_2n +c$ shuffles, that is when $a=2^c n$,
\refprop{TV-bound} shows that the total variation distance is
approximately (with $C=2^c$)
\[ 
C\log \Big( C (e^{1/C}-1)\Big) + \frac{1-C \log (e^{1/C}-1)}{(e^{1/C}-1)}.
\]
Thus when $c$ is `large and negative,' the total variation is close to
$1$, and when $c$ is large and positive, the total variation is close
to $0$.  Thus total variation and separation converge at the same
rate. This is an asymptotic result and, for example, \reftab{AS}
supports this.

Similar, but more demanding, calculations show that if the ace of
spades starts at position $i$, and $\max(i/n, (n-i)/n) \geq A > 0$ for
some fixed positive $A$, then $\frac{1}{2} \log_2 n$ shuffles suffice
for convergence in any of the metrics. We omit further details.

%%%%%%%%%%%%%%%%%%%%%%%%%%%%%%%%%%%%%%%%%%%%%%%%%%%%%%%%%%%%
\section{Separation distance for the general case}
%%%%%%%%%%%%%%%%%%%%%%%%%%%%%%%%%%%%%%%%%%%%%%%%%%%%%%%%%%%%
\label{sec:deck}

A main result of Bayer and Diaconis \cite{BaDi1992} is the simple
formula for an $a$-shuffle of a deck of $n$ distinct cards:
\begin{equation}
  Q_a (\sigma) = \frac{1}{a^n}\binom{n+a-r}{n},
\label{eqn:BD}
\end{equation}
where $r = r(\sigma)$ is the number of rising sequences in $\sigma$,
equivalently one more than the number of descents in
$\sigma^{-1}$. This formula allows simple closed form expressions for
a variety of distances as well as asymptotic analysis.

In this section we work with general decks containing $D_i$ cards
labelled $i$, $1 \leq i \leq m$. The formulae of this section hardly
resemble the elegant expression above. Further, we only give precise
formula for the least likely deck. The following lemma shows that this
deck, where the separation distance is achieved, is the reverse the
initial deck configuration. This is equivalent to Theorem 2.1
from \cite{CoVi2006}.

\begin{proposition}
  Let $D$ be a deck as above. After an $a$-shuffle of the deck with
  $1$'s on top down to $m$'s on bottom, the least likely configuration
  is the reverse deck $w^*$ with $m$'s on top down to $1$'s on the
  bottom. 
\label{prop:sepw}
\end{proposition}

\begin{proof}
  The only cuts of the initial deck resulting in $w^*$ are those
  containing no pile with distinct letters. For all such cuts, each
  rearrangement of the deck is equally likely to occur. 
\end{proof}

While finding a completely general formula for $Q_a(w)$ for arbitrary
$w$ is infeasible, below we do this for $w^*$.

\begin{theorem}
  Consider a deck with $n$ cards and $D_i$ cards labeled $i$,
  $i=1,\ldots,m$. Then the separation distance after an $a$-shuffle of
  the sorted deck ($1$'s followed by $2$'s, etc) is given by
  \begin{displaymath}
    \SEP(a) = 
    1 \!-\! \frac{1}{a^n} \binom{n}{D_1\ldots D_m}
    \hspace{-2.5em} \sum_{\rule{0pt}{3ex} 0=k_0 < \cdots < k_{m-1} < a}
    \hspace{-2em} (a\!-\!k_{m\!-\!1})^{D_m} \! \prod_{j=1}^{m-1} 
    \left( (k_j \!-\! k_{j\!-\!1})^{D_j} \!-\! (k_j \!-\! k_{j\!-\!1}
      \!-\! 1)^{D_j} \right) .
  \end{displaymath}
\label{thm:sepG}
\end{theorem}

\begin{proof}
  From the analysis in the proof of \refprop{sepw}, $Q_a(w^*)$ is
  given by
  \begin{displaymath}
    Q_a(w^*) = 
    \sum_{\substack{A_1 + \cdots + A_a = n \\ A \ \mathrm{refines} \ D}}
    \frac{1}{a^n} \binom{n}{A_1,\ldots,A_a}
    \frac{1}{\binom{n}{D_1,\ldots,D_m}},
  \end{displaymath}
  where `$A$ refines $D$' means there exist indices
  $k_1,\ldots,k_{m-1}$ such that $A_1+\cdots+A_{k_1} = D_1$ and, for
  $i=2,\ldots,m-1$, $A_{k_{i-1} + 1} + \cdots + A_{k_i} = D_i$. Taking
  the $k_i$'s to be minimal, the expression for $Q_a(w^*)$ simplifies
  to
  \begin{equation}
    \frac{1}{a^n} \hspace{-1ex} \sum_{0=k_0 < \cdots <
      k_{m-1} < a} \hspace{-2em} (a\!-\!k_{m-1})^{D_m} \prod_{j=1}^{m-1}
    \left( (k_j \!-\! k_{j-1})^{D_j} - (k_j \!-\! k_{j-1} \!-\! 1)^{D_j} \right).
  \label{eqn:sepG}
  \end{equation}
  The result now follows from \refprop{sepw}.
\end{proof}

%%The Bayer-Diaconis formula (\ref{eqn:BD}) for the reverse deck yields
%%the following formula for separation distance when all $n$ cards are
%%distinct,
%%\begin{equation}
%%  \SEP(a) = 1 - \frac{n!}{a^n} \binom{a}{n}.
%%\end{equation}
%%Comparing this with the formula in \refthm{sepG} is somewhat
%%disheartening. Nevertheless, this formula allows us to make some
%%calculations summarized in \reftab{sep}.

\begin{table}[ht]
  \begin{center}
    \caption{\label{tab:sep}Separation distance for $k$ shuffles of $52$ cards.}
    \begin{tabular}{c | cccc cccc cccc}
      \hline
      k & 1 & 2 & 3 & 4 & 5 & 6 & 7 & 8 & 9 & 10 & 11 & 12 \\ \hline
      \\ [-.5\vsp] 
      BD-92 & 1.00 & 1.00 & 1.00 & 1.00 & 1.00 & 1.00 &
      1.00 & .995 & .928 & .729 & .478 & .278 \\ [.5\vsp]
      blackjack & 1.00 & 1.00 & 1.00 & 1.00 & .999 & .970 &  &
       &  &  &  & \\ [.5\vsp]  
      $\clubsuit{\red\diamondsuit\heartsuit}\spadesuit$ & 1.00
      & .997 & .997 & .976 & .884 & .683 & .447 & .260 & .140 & .073
      &  &  \\ [.5\vsp]  
      A$\spadesuit$ & 1.00 & 1.00 & .993 & .875 & .605 & .353 & .190 &
      .098 & .050 & .025 & .013 & .006 \\ [.5\vsp] 
      {\red red}black & .890 & .890 & .849 & .708 & .508 & .317 & .179 &
      .095 & .049 & .025 & .013 & .006 \\ [.3\vsp]
      \raisebox{-.5ex}{\makebox[0pt]{\includegraphics[height=3.5ex]{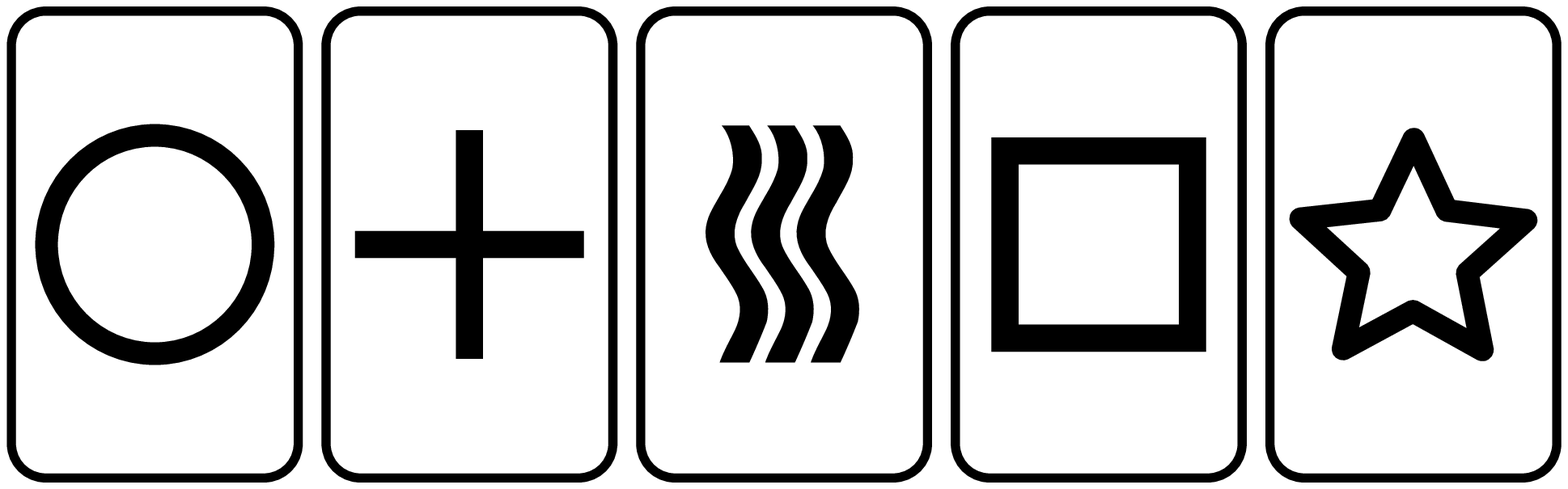}}}
      & 1.00 & 1.00 & .993 & .943 & .778 & .536 & .321 & .177
      &  &  &  & \\
      \hline
    \end{tabular}
  \end{center}
\end{table}

\vspace{\baselineskip}
\begin{rmkT}{sep}
  We calculate $\SEP$ after repeated $2$-shuffles for various decks
  using \refthm{sepG}: (blackjack) $9$ ranks with $4$ cards each and
  another rank with $16$ cards; ($\clubsuit \diamondsuit \heartsuit
  \spadesuit$) $4$ distinct suits of $13$ cards each; (A$\spadesuit$)
  the ace of spades and $51$ other cards; (redblack) a two color deck
  with $26$ of either color; and
  (\raisebox{-.8ex}{\includegraphics[height=3ex]{zener.eps}}) a deck
  with 5 cards in each of 5 suits. The missing entries in \reftab{sep}
  highlight the limitations of exact calculations using \refthm{sepG}.
\end{rmkT}

\begin{remark}
  Comparing the data in \reftab{sep} for $\mathrm{A}\spadesuit$ and
  $\mathrm{redblack}$ shows that these two cases are remarkably
  similar. Indeed, both cases exhibit the same asymptotic behavior,
  which is remarkable since the $\mathrm{A}\spadesuit$ has a state
  space of size $52$ while $\mathrm{redblack}$ has a state space of
  size around $5 \times 10^{14}$.
\end{remark}

Now we derive a basic asymptotic tool which allows asymptotic
approximations for general decks. 

%%Again, full details may be found in \cite{ADS2008}. As motivation,
%%consider again the case of one card mixing, i.e. begin with $n$ cards
%%with the ace of spaces at the bottom of the initial deck. Recall from
%%\refcor{Qai} the chance that the ace of spades is at position $i$ from
%%the top after an $a$-shuffle. By \refprop{sepw}, we have
%%\begin{equation}
%%  \SEP(a) = 1 - nQ_a(1) = 1 - \frac{n}{a^n} \sum_{k=1}^{a}(k-1)^{n-1},
%%  \label{eqn:sepi}
%%\end{equation}
%%with the convention $0^0=1$. Exact calculations when $n=52$ are given
%%in \reftab{sep}. 

\begin{proposition} Let $m\ge 2$ and $a$ be natural numbers, let
  $\xi_1$, $\ldots$, $\xi_m$ be real numbers in $[0,1]$.  Let $r_1$,
  $\ldots$, $r_m$ be natural numbers all at least $r\ge 2$.  Let
  \begin{displaymath}
    S_m(a; \underline{\xi}, \underline{r}) =
    \sum_{\substack{a_1,\ldots,a_m \ge 0 \\ a_1+\ldots+a_m=a}}
    (a_1+\xi_1)^{r_1} \cdots (a_m+\xi_m)^{r_m}.
  \end{displaymath}
  Then 
  \begin{displaymath}
    \begin{array}{l}
      \displaystyle{\Big| S_{m}(a;\underline{\xi},\underline{r}) 
        -\frac{r_1! \cdots r_m!}{(r_1+\ldots+r_m+m-1)!}
        (a+\xi_1+\ldots+\xi_m)^{r_1+\ldots+r_m+m-1}\Big|} \\
      \displaystyle{\leq r_1! \cdots r_m! \sum_{j=1}^{m-1}
        \binom{m-1}{j} \Big(\frac{1}{3 (r-1)}\Big)^j  \frac{ (a
          +\xi_1+\ldots+\xi_m)^{r_1+\ldots+r_m
            +m-1-2j}}{(r_1+\ldots+r_m+m-1-2j)!}}.
    \end{array}
  \end{displaymath}
  \label{prop:GPS}
\end{proposition} 
 
Consider a general deck of $n$ cards with $D_i$ cards labelled $i$. We
use \refprop{GPS} to find asymptotics for the separation distance
given in \refthm{sepG}. The following is our `rule of thumb.'
 
\begin{theorem} 
  For a deck of $n$ cards as above, suppose $D_i\ge d \ge 3$ for all
  $1\le i \le m$. Then we have
  \[ \SEP(a) = 
  1- (1+\eta) \frac{a^{m-1}}{(n+1)\cdots(n+m-1)} \sum_{j=0}^{m-1} (-1)^j
  \binom{m-1}{j} \Big(1-\frac ja\Big)^{n+m-1} ,
  \]
  where $\eta$ is a real number satisfying 
  \[ 
  |\eta| \le \Big(1 +\frac{n^2}{3(d-2)(a-m+1)^2}\Big)^{m-1} -1.
  \]
  \label{thm:ROT}  
\end{theorem} 

\begin{proof}   
  To evaluate the expression in \refthm{sepG}, we require an
  understanding of
  \begin{displaymath}
      \sum_{\substack{a_1+ \ldots + a_m =a \\ a_j \ge 1}} \!\!\!
      a_m^{D_m} \prod_{j=1}^{m-1} (a_j^{D_j}-(a_j-1)^{D_j})
      = \int_0^1 \cdots \int_0^1 \sum_{\substack{a_1
          +\ldots+a_m=a \\ a_j\ge 1}} \!\!\! a_m^{D_m} \prod_{j=1}^{m-1}
      \Big( D_j (a_j-1+\xi_j)^{D_j-1} d\xi_j \Big).
  \end{displaymath}    
  We now invoke \refprop{GPS}.  Thus the above equals for some
  $|\theta| \le 1$
  \begin{displaymath}
    \begin{array}{l}
      \displaystyle{%
        \prod_{j=1}^{m} D_j!  \int_0^1 \cdots \int_0^1 \Big( 
        \frac{(a-(m-1)+\xi_1+\ldots+\xi_{m-1})^{n}}{n!} + }\\
      \displaystyle{+ \theta \sum_{j=1}^{m-1} \!\!
        \binom{m\!-\!1}{j} \! \Big(\frac{1}{3(d\!-\!2)}\Big)^j
        \frac{(a\!-\!(m\!-\!1)\!+\!\xi_1\!+\!\ldots\!+\!
          \xi_{m-1})^{n-2j}}{(n-2j)!}\Big) d\xi_1 \!\cdots\! d\xi_{m-1}.}
    \end{array}
  \end{displaymath}
  We may simplify the above as
  \begin{displaymath} 
    \Big( 1+ \theta \Big\{
    \Big(1+\frac{n^2}{3(d\!-\!2)(a\!-\!m\!+\!1)^2}\Big)^{m\!-\!1}
    \!\!\! -1\Big\} \Big) \frac{D_1! \cdots D_m!}{n!}
    \int_0^1\!\!\!\ldots\!\!\int_0^1 (a\!-\!m\!+\!1 +\xi_1+\cdots+\xi_{m-1})^{n} 
    d\xi_1 \cdots d\xi_{m-1},
  \end{displaymath}
  and evaluating the integrals above this is 
  \begin{displaymath}
    \left( 1+ \theta \Big\{\left( 1+\frac{n^2}{3(d-2)(a-m+1)^2}
      \right)^{m-1}-1\Big\} \right)
    \frac{D_1! \cdots D_m!}{n!} \sum_{j=0}^{m-1} (-1)^j
    \binom{m-1}{j} (a-j)^{n-m+1}.
  \end{displaymath}
  The Theorem follows.
\end{proof}

For simplicity we have restricted ourselves to the case when each pile
has at least three cards.  With more effort we could extend the
analysis to include doubleton piles.  The case of some singleton piles
needs some modifications to our formula, but this variant can also be
worked out. Below we use our rule of thumb to calculate separation for
the same decks as in \reftab{sep}.

\begin{table}[ht]
  \begin{center}
    \caption{\label{tab:thumb} Rule of Thumb for the separation
      distance for $k$ shuffles of $52$ cards.}
    \begin{tabular}{c | cccc cccc cccc}
      \hline
     k & 1 & 2 & 3 & 4 & 5 & 6 & 7 & 8 & 9 & 10 & 11 & 12 \\ \hline
     \\ [-.5\vsp]
     BD-92 & 1.00 & 1.00 & 1.00 & 1.00 & 1.00 & 1.00 & 1.00 &
     .995 & .928 & .729 & .478 & .278 \\ [.5\vsp]
     blackjack & 1.00 & 1.00 & 1.00 & 1.00 & .999 & .970 & .834 &
     .596 & .366 & .204 & .108 & .056 \\ [.5\vsp]  
     $\clubsuit{\red\diamondsuit\heartsuit}\spadesuit$ & 1.00
     & 1.00 & .997 & .976 & .884 & .683 & .447 & .260 & .140 & .073
     & .037 & .019 \\ [.5\vsp]  
     {\red red}black & .962 & .925 & .849 & .708 & .508 & .317 & .179 &
     .095 & .049 & .025 & .013 & .006 \\%[.5\vsp]
     \raisebox{-.5ex}{\makebox[0pt]{\includegraphics[height=3.3ex]{zener.eps}}}
     &  1.00 & 1.00 & .993 & .943 & .778 & .536 & .321 & .177 & .093 &
     .048 & .024 & .012 \\ \hline 
   \end{tabular}
 \end{center}
\end{table}

\begin{rmkT}{thumb}
  The first row gives exact results from the Bayer-Diaconis formula
  for the full permutation group. The other numbers are from the rule
  of thumb. Roughly, the single card or red-black numbers suggest that
  half the usual number of shuffles suffice. The Black-Jack
  (equivalently Baccarat) numbers suggest a savings of two or three
  shuffles, and the suit numbers lie in between. The final row is the
  rule of thumb for the Zener deck with 25 cards, 5 cards for each of
  5 suits.
\end{rmkT}

While asymptotic, \refthm{ROT} is astonishingly accurate for decks of
practical interest. For instance, comparing exact calculations in
\reftab{sep} with approximations using this rule of thumb in
\reftab{thumb} shows that after only $3$ shuffles, the numbers agree
to the given precision. Moreover, the simplicity of the formula in
\refthm{ROT} allows much further computations than are possible using
the formula in \refthm{sepG}.

We now give a heuristic for why our rule of thumb is numerically so
accurate.  For $k\geq 0$, define
\begin{displaymath}
  f_k(z) = \sum_{r=0}^{\infty} r^{k} z^r = \frac{A_k(z)}{(1-z)^{k+1}},
\end{displaymath}
where $A_k(z)$ denotes the $k$-th Eulerian polynomial. The sum over
$a_1$, $\ldots$, $a_m$ appearing in our proof of \refthm{ROT} is
simply the coefficient of $z^{a}$ in the generating function
$(1-z)^{m-1}f_{D_1}(z)\cdots f_{D_m}(z)$.  Our rule of thumb may be
interpreted as saying that
\begin{equation}
(1-z)^{m-1} f_{D_1}(z)\cdots f_{D_m}(z) \approx 
\frac{D_1! \cdots D_{m}!}{(n+m-1)!} (1-z)^{m-1} f_{n+m-1}(z).
\label{eqn:genROT}
\end{equation}
To explain the sense in which \refeq{genROT} holds, note that $f_k(z)$
extends meromorphically to the complex plane, and it has a pole of
order $k+1$ at $z=1$.  Moreover it is easy to see that $f_k(z) -
k!/(1-z)^{k+1}$ has a pole of order at most $k$ at $z=1$.  Therefore,
the LHS and RHS of \refeq{genROT} have poles of order $n+1$ at $z=1$,
and their leading order contributions match.  Therefore the difference
between the RHS and LHS of \refeq{genROT} has a pole of order at most
$n$ at $z=1$.  But in fact, this difference can have a pole of order
at most $n-d$ at $z=1$, and thus the approximation in \refeq{genROT}
is tighter than what may be expected {\sl a priori}.  To obtain our
result on the order of the pole, we record that one can show
\begin{displaymath}
  f_k(z) = \frac{k!}{(1-z)^{k+1}} \Big( \frac{(z-1)}{ \log z}
  \Big)^{k+1} +\zeta(-k) + O(1-z). 
\end{displaymath}

%%%%%%%%%%%%%%%%%%%%%%%%%%%%%%%%%%%%%%%%%%%%%%%%%%%%%%%%%%%%
\section{Gilbreath principle at work}
%%%%%%%%%%%%%%%%%%%%%%%%%%%%%%%%%%%%%%%%%%%%%%%%%%%%%%%%%%%%
\label{sec:initial}

Conger and Viswanath note that the initial configuration can affect
the speed of convergence to stationary. Perhaps this is most striking
in the case of \refsec{Aspades} where a single card is tracked. Recall
\reftab{AS}, giving calculations for the distinguished card beginning
at the bottom of a deck of $52$ cards. In contrast, \reftab{AD} gives
calculations for the distinguished card starting in the middle, at
position $26$. For the latter, both total variation and separation are
indistinguishable from zero after only four shuffles.

\begin{table}[ht]
  \begin{center}
    \caption{\label{tab:AD}Distance to uniformity for a single card
      starting at the middle of a $52$ card deck.}
    \begin{tabular}{c | rrrr}
      \hline
       & 1 & 2 & 3 & 4 \\ \hline \\
       [-.5\vsp]  
      $TV$ & .494 & .152 & .001 & .000 \\ [.5\vsp]
      $\SEP$ & 1.00 & .487 & .003 & .000 \\ \hline
    \end{tabular}
  \end{center}
\end{table}

Consider next a deck with $n$ red and $n$ black cards. First take the
starting condition of all reds atop all blacks. If the initial cut is
at $n$ (the most likely value) then the red-black pattern is perfectly
mixed after a single shuffle. More generally, the chance of the deck
$w$ resulting from a single $2$-shuffle of a deck with $n$ red cards
atop $n$ black cards is given by
\begin{displaymath}
    Q_2(w) = \frac{1}{2^{2n}} \left( 2^{\head(w)} + 2^{\tail(w)} - 1 \right),
\end{displaymath}
where $\head(w)$ is the number of red cards before the first black
card and $\tail(w)$ is the number of black cards after the final red
card; see \cite{ADS2008}. In particular, the total variation after a
single $2$-shuffle is 
\begin{equation}
  \left\|Q_2 - U\right\|_{TV} = \frac{1}{2} \left(
    \left(\frac{2^{n+1}\!-\!1}{2^{2n}} - \frac{1}{\binom{2n}{n}}\right)
    + \sum_{i=0}^{n-1} \sum_{j=0}^{n-1} \left| \frac{2^i \!+\! 2^j \!-\!
        1}{2^{2n}} - \frac{1}{\binom{2n}{n}} \right| \binom{2n\!-
      (i\!+\!j\!+\!2)}{n\!-(i\!+\!1)} \right)
\label{eqn:2sh}
\end{equation}
Evaluating this formula for $2n=52$ give a total variation of $0.579$.

Now take the starting condition to alternate red black red black,
etc. As motivation, we recall a popular card trick: Begin with a deck
of $2n$ cards arranged alternately red, black, red, black, etc. The
deck may be cut any number of times. Have the deck turned face up and
cut (with cuts completed) until one of the cuts results in the two
piles having cards of opposite color uppermost. At this point, ask one
of the participants to riffle shuffle the two piles together. The
resulting arrangement has the top two cards containing one red and one
black, the next two cards containing one red and one black, and so on
throughout the deck. This trick is called the Gilbreath Principle
after its inventor, the mathematician Norman Gilbreath. It is
developed, with many variations, in Chapter 4 of \cite{Gardner1966}.
From the trick we see that beginning with an alternating deck severely
limits the possibilities. Analyzing the trick reveals the following
formula,
\begin{equation}
  2^{2n} \cdot Q_2(w) = \left\{ \begin{array}{cl}
      2^{n-1} + 2^n & \mbox{if $w$ is the initial alternating deck}, \\
      2^{n-1}       & \mbox{if $w$ can result from an odd cut}, \\
      2^n          & \mbox{if $w$ can result from an even cut}, \\
      0            & \mbox{otherwise},
    \end{array} \right.
  \label{eqn:Q2-2}
\end{equation}
where an odd (resp. even) cut refers to the parity of cards in either
pile. From this we compute
\begin{equation}
  \left\|Q_2 - U\right\|_{TV} = \frac{1}{2} \left( 1 - \frac{2^n +
      2^{n-1} - 1}{\binom{2n}{n}} \right),
  \label{eqn:TValt}
\end{equation}
which goes to $.5$ exponentially fast as $n$ goes to infinity, and
indeed is already $.500$ for $2n=52$. In contrast, starting with reds
above blacks, asymptotic analysis of \refeq{2sh} shows that the total
variation tends to $1$ after a single shuffle when $n$ is large. Thus
again an alternating start leads to faster mixing.

%%%%%%%%%%%%%%%%%%%%%%%%%%%%%%%%%%%%%%%%%%%%%%%%%%%%%%%%%%%%
  \bibliographystyle{abbrv} 
  \bibliography{../../references}

\begin{thebibliography}{10}

\bibitem{AlDi1986}
D.~Aldous and P.~Diaconis.
\newblock Shuffling cards and stopping times.
\newblock {\em Amer. Math. Monthly}, 93(5):333--348, 1986.

\bibitem{ADS2008}
S.~Assaf, P.~Diaconis, and K.~Soundararajan.
\newblock A rule of thumb for riffle shuffling.
\newblock Preprint, 2008.

\bibitem{BaDi1992}
D.~Bayer and P.~Diaconis.
\newblock Trailing the dovetail shuffle to its lair.
\newblock {\em Ann. Appl. Probab.}, 2(2):294--313, 1992.

\bibitem{Ciucu1998}
M.~Ciucu.
\newblock No-feedback card guessing for dovetail shuffles.
\newblock {\em Ann. Appl. Probab.}, 8(4):1251--1269, 1998.

\bibitem{CoVi2006}
M.~Conger and D.~Viswanath.
\newblock Riffle shuffles of decks with repeated cards.
\newblock {\em Ann. Probab.}, 34(2):804--819, 2006.

\bibitem{CoVi}
M.~Conger and D.~Viswanath.
\newblock Shuffling cards for blackjack, bridge, and other card games.
\newblock Preprint, 2006.

\bibitem{CoVi2007}
M.~Conger and D.~Viswanath.
\newblock Normal approximations for descents and inversions of permutations of
  multisets.
\newblock {\em J. Theoret. Probab.}, 20(2):309--325, 2007.

\bibitem{Diaconis2003}
P.~Diaconis.
\newblock Mathematical developments from the analysis of riffle shuffling.
\newblock In {\em Groups, combinatorics \& geometry (Durham, 2001)}, pages
  73--97. World Sci. Publ., River Edge, NJ, 2003.

\bibitem{DiFu2008}
P.~Diaconis and J.~Fulman.
\newblock Carries, shuffling and an amazing matrix.
\newblock Preprint, 2008.

\bibitem{Gardner1966}
M.~Gardner.
\newblock {\em {M}artin {G}ardner's New Mathematical Diversions from
  {S}cientific {A}merican}.
\newblock Simon \& Schuster, New York, 1966.

\bibitem{Gilbert1955}
E.~Gilbert.
\newblock Theory of shuffling.
\newblock Technical memorandum, Bell Laboratories, 1955.

\bibitem{Holte1997}
J.~M. Holte.
\newblock Carries, combinatorics, and an amazing matrix.
\newblock {\em Amer. Math. Monthly}, 104(2):138--149, 1997.

\bibitem{Reeds1976}
J.~Reeds.
\newblock Theory of shuffling.
\newblock Unpublished manuscript, 1976.

\end{thebibliography}
%%%%%%%%%%%%%%%%%%%%%%%%%%%%%%%%%%%%%%%%%%%%%%%%%%%%%%%%%%%%

\end{document}